\documentclass[11pt]{article}
\usepackage{amsmath,amsfonts,amssymb,amsthm,graphicx,verbatim,subfig,mathtools,wasysym}
\usepackage[all]{xy}
\usepackage{dsfont,enumitem}
\usepackage[a4paper,left=25mm]{geometry}
\ifx\pdftexversion\undefined

\usepackage[a4paper,colorlinks,linkcolor=black,filecolor=black,citecolor=black,urlcolor=black,pdfstartview=FitH]{hyperref}
\fi

\font\sixbb=msbm6
\font\eightbb=msbm8
\font\twelvebb=msbm10 scaled 1095
\newfam\bbfam
\textfont\bbfam=\twelvebb \scriptfont\bbfam=\eightbb
                           \scriptscriptfont\bbfam=\sixbb
\def\bb{\fam\bbfam\twelvebb}
\newcommand{\Rea}{\mathbb{R}}
\newcommand{\Com}{\mathbb{C}}

\newcommand{\FF}{{\bb F}}

\newcommand{\tr}{\text{tr}}
\newcommand{\cll}{\mathcal{L}}
\newcommand{\crr}{\mathcal{R}}

\newtheorem{theorem}{\bf Theorem}[section]
\newtheorem{claim}[theorem]{\bf Claim}
\newtheorem{conjecture}[theorem]{\bf Conjecture}
\newtheorem{proposition}[theorem]{\bf Proposition}

\newtheorem{definition}[theorem]{\bf Definition}
\newtheorem{remark}[theorem]{\bf Remark}

\newcommand{\enp}{\begin{flushright} $\Box$ \end{flushright}}
\newcommand{\beq}[0]{\begin{equation}}
\newcommand{\enq}[0]{\end{equation}}

\newcommand{\thh}{\tilde{H}}
\newcommand{\supp}{{\rm supp}}
\newcommand{\cf}{{\cal F}}
\newcommand{\pr}{{\rm Pr}}

\newcommand{\lk}{{\rm lk}}

\newcommand{\img}{{\rm Im} \,}

\newcommand{\one}{\mathds{1}}

\newcommand{\phat}{\widehat{\phi}}

\newcommand{\ghat}{\widehat{G}}

\newcommand{\ccs}{\mathcal{S}}

\newcommand{\mb}[1]{{\mathbf #1}}

\newcommand{\rmc}{\mathcal{C}}

\newcommand{\prob}{{\rm Pr}}
\newcommand{\ddh}{\Delta_{h-1}}
\newcommand{\csy}{{\rm csy}}
\newcommand{\sph}{\textsc{S}}
\newcommand{\ueps}{\underline{\epsilon}}

\title{Random Balanced Cayley Complexes}

\author{Roy Meshulam\thanks{Department of Mathematics,
Technion, Haifa 32000, Israel. e-mail:
meshulam@technion.ac.il~. Supported by ISF grant 686/20.}}
\begin{document}
\maketitle
\begin{abstract}
Let $G$ be a finite group of order $n$ and for $1 \leq i \leq k+1$ let $V_i=\{i\} \times G$.
Viewing each $V_i$ as a $0$-dimensional complex, let $Y_{G,k}$ denote the simplicial join $V_1*\cdots*V_{k+1}$.
For $A \subset G$ let $Y_{A,k}$ be the subcomplex of $Y_{G,k}$ that contains the $(k-1)$-skeleton
of $Y_{G,k}$ and whose $k$-simplices are
all $\{(1,x_1),\ldots,(k+1,x_{k+1})\} \in Y_{G,k}$ such that $x_1\cdots x_{k+1} \in A$.
Let $L_{k-1}$ denote the reduced $(k-1)$-th Laplacian
of $Y_{A,k}$, acting on the space $C^{k-1}(Y_{A,k})$ of real valued $(k-1)$-cochains of $Y_{A,k}$.  The $(k-1)$-th spectral gap $\mu_{k-1}(Y_{A,k})$ of $Y_{A,k}$ is the minimal eigenvalue of $L_{k-1}$.
\\
The following $k$-dimensional analogue of the Alon-Roichman theorem is proved:
Let $k \geq 1$ and $\epsilon>0$ be fixed and let $A$ be a random subset of $G$ of size $m=
\left\lceil\frac{10 k^2\log D}{\epsilon^2}\right\rceil$ where $D$ is the sum of the degrees of the complex irreducible representations of $G$.
Then
\[
\pr\big[~\mu_{k-1}(Y_{A,k}) < (1-\epsilon)m~\big] =O\left(\frac{1}{n}\right).
\]
\end{abstract}
\section{Introduction}
\label{s:intro}

The Laplacian $L(\rmc)$ of a graph
$\rmc=(V,E)$ is the  $V \times V$ positive semidefinite matrix whose $(u,v)$ entry is given by
\begin{equation*}
L(\rmc)_{uv}= \left\{
\begin{array}{cl}
\deg_C(u) & u=v, \\
-1 & \{u,v\} \in E, \\
0 & {\rm otherwise}.
\end{array}
\right.
\end{equation*}
Let $0=\lambda_1(\rmc) \leq \lambda_2(\rmc) \leq \cdots \leq \lambda_{|V|}(\rmc)$ be the eigenvalues of $L(\rmc)$.
The second smallest eigenvalue $\lambda_2(\rmc)$, called the \emph{spectral gap} of $\rmc$,
is a parameter of central importance in a variety of problems. In
particular it controls the expansion properties of $\rmc$ and the
convergence rate of a random walk on $\rmc$ (see e.g. chapters XIII and IX in \cite{Bollobas98}).

Throughout the paper, let $G$ denote a finite group of order $n$
and let $\ghat=\{\rho\}$ be the set of irreducible
unitary representations of $G$, where $\rho:G \rightarrow U(d_{\rho})$. Let $D(G)=\sum_{\rho \in \ghat} d_{\rho}$.
Let $\one \in \ghat$ denote the trivial representation of $G$ and let $\ghat_+=\ghat\setminus \{\one\}$.

Let $T \subset G$ be symmetric subset, i.e. $T=T^{-1}$.
The \emph{Cayley graph} $\rmc(G,T)$ of $G$ with respect to $T$
is the graph on the vertex set $G$ with edge set $\left\{\{g,gt\}: g \in G, t \in T \right\}$.
The seminal  Alon-Roichman theorem \cite{AR94} is concerned with the
expansion of Cayley graphs with respect to random sets of generators.
\begin{theorem}[Alon-Roichman]
\label{t:AR}
For any $\epsilon>0$ there exists a constant $c(\epsilon)>0$ such that for any group $G$,
if $S$ is a random subset of $G$ of size
$\lceil c(\epsilon) \log |G| \rceil$ and $m=|S \cup S^{-1}|$, then $\lambda_2(\rmc\left(G,S \cup S^{-1})\right)$ is
asymptotically almost surely (a.a.s.) at least $(1-\epsilon)m$.
\end{theorem}
\noindent
\begin{remark}
Landau and Russell \cite{Landau-Russell04} and independently Loh and Schulman \cite{Loh-Schulman04} have obtained an improvement on Theorem \ref{t:AR} by showing that the $\log |G|$ factor in the bound on $|S|$ can be replaced by $\log D(G)$. While this does not change the logarithmic dependence of $|S|$ on $|G|$,
it does often lead to an improvement of the constant $c(\epsilon)$.
\end{remark}

This paper is concerned with higher dimensional counterparts of Theorem \ref{t:AR}.
We briefly recall the relevant terminology (see section \ref{s:hodge} for
details). For a simplicial complex $X$ and $k \geq -1$
let $X^{(k)}$ denote the $k$-dimensional skeleton of $X$.
For $k \geq -1$ let $C^k(X)$ denote the space of real valued
simplicial $k$-cochains of $X$ and let $d_k:C^k(X)
\rightarrow C^{k+1}(X)$ denote the coboundary operator. For $k
\geq 0$ define the reduced $k$-th Laplacian of $X$ by
$L_k(X)=d_{k-1}d_{k-1}^*+d_k^*d_k$.  The minimal eigenvalue of $L_k(X)$, denoted by $\mu_k(X)$, is the \emph{$k$-th spectral gap} of $X$.

The following $k$-dimensional abelian version of Theorem \ref{t:AR} was obtained in \cite{BAM20}.
Let $H$ be an additively written abelian group of order $h$ and let $k \leq h$. Let $\ddh$ denote the $(h-1)$-simplex on the vertex set $H$. The \emph{Sum Complex} $X_{A,k}$ associated with a subset $A \subset H$ is  the $k$-dimensional simplicial complex obtained
by taking the full $(k-1)$-skeleton of $\ddh$ together with all $(k+1)$-subsets $\sigma \subset H$ that satisfy $\sum_{x \in \sigma} x \in A$.
\begin{theorem}[\cite{BAM20}]
\label{t:bam20}
Let $k \geq 1$ and $\epsilon>0$ be fixed and let $A$ be a random subset of $H$ of size $m=\left\lceil\frac{4k^2\log h}{\epsilon^2}\right\rceil$. Then
\begin{equation*}
\pr\big[~\mu_{k-1}(X_{A,k}) < (1-\epsilon)m~\big] =O\left(\frac{1}{n}\right).
\end{equation*}
\end{theorem}
\begin{remark}
See \cite{LMR10,M14} for more on sum complexes and their cohomology, .
\end{remark}
In the present paper we study a different model of Cayley complexes associated with
subsets of general finite groups and
obtain a new high dimensional analogue of Theorem \ref{t:AR}.
Recall that $G$ is a finite group of order $n$ and let $k \geq 1$. For $1 \leq i \leq k+1$ let
$V_i=\{i\} \times G$.  Let $Y_{G,k}$ denote the simplicial join $V_1* \cdots * V_{k+1}$,
where each $V_i$ is viewed as $0$-dimensional complex. Thus $Y_{G,k}$ is homotopy equivalent to an
$N$-fold wedge $\bigvee^N \sph^k$ of $k$-dimensional spheres, where $N=(n-1)^{k+1}$.
For $\emptyset \neq A \subset G$ let
\begin{equation*}
P_{A,k}=\{\mb{x}=(x_0,\ldots,x_k) \in G^{k+1}: x_0\cdots x_{k} \in A\}.
\end{equation*}
\noindent
The \emph{balanced $k$-dimensional Cayley Complex} associated with
$A$ is the simplicial
complex $Y_{G,k}^{(k-1)} \subset Y_{A,k} \subset Y_{G,k}$ whose $k$-dimensional simplices are
$\left\{ (1,y_1),\ldots,(k+1,y_{k+1})\right\}$ where $(y_1,\ldots,y_{k+1}) \in P_{A,k}$.
\noindent
Let $1_A$ denote the indicator function of $A \subset G$, i.e. $1_A(x)=1$ if $x \in A$ and $1_A(x)=0$ otherwise. For a representation $\rho:G \rightarrow GL_d(\Com)$ let
$\widehat{1_A}(\rho)=\sum_{x \in A} \rho(x) \in M_d(\Com)$ be the Fourier transform of $1_A$ at $\rho$ (see
section \ref{s:fourier} for details).
For a matrix $T \in M_d(\Com)$ let $\|T\|=\max_{\|v\|=1} \|Tv\|$ denote the spectral norm of $T$.
Let
%\begin{equation*}
$\nu(A)=\max_{\rho \in \ghat_+} \|\widehat{1_A}(\rho)\|$.
%\end{equation*}
Our first result is a lower bound on the $(k-1)$-th spectral gap of $Y_{A,k}$ in terms of $\nu(A)$.
\begin{theorem}
\label{t:specgap}
\begin{equation*}
\label{e:specgap}
\mu_{k-1}(Y_{A,k}) \geq |A|-k \cdot \nu(A).
\end{equation*}
\end{theorem}
Our main result is the following $k$-dimensional analogue of the Alon-Roichman Theorem.
\begin{theorem}
\label{t:logD}
Let $k$ and $\epsilon>0$ be fixed. Suppose that $|G|=n > 10^6 \left(\frac{k}{\epsilon}\right)^8$ and let $A$ be a random subset of $G$ of size
$m=\left\lceil\frac{9k^2\log D(G)}{\epsilon^2}\right\rceil$. Then
\begin{equation*}
\label{e:sprob}
\pr\big[~\mu_{k-1}(Y_{A,k}) < (1-\epsilon)m~\big] < \frac{6}{n}.
\end{equation*}
\end{theorem}
\begin{remark}
While there are some similarities between sum complexes and balanced Cayley complexes, the analysis of $Y_{A,k}$ in the present paper requires some additional ideas, including the use of the non-abelian Fourier transform and of Garland's eigenvalue estimates \cite{Garland73}.
\end{remark}
The paper is organized as follows. In Section \ref{s:hodge} we review
some basic properties of high dimensional Laplacians and their eigenvalues, including Garland's
lower bound for the higher spectral gaps.
In Section \ref{s:sykg} we compute the spectra of various Laplacians of the skeleta of $Y_{n,k}$ and deduce
a variational characterization (Proposition \ref{p:sgapy}) of $\mu_{k-1}(Y)$
for subcomplexes $Y_{G,k}^{(k-1)} \subset Y \subset Y_{G,k}$. In Section \ref{s:fourier} we briefly recall
the definition and some basic properties of the Fourier transform on finite groups. In Section \ref{s:sgyak} we prove Theorem
\ref{t:specgap}. This bound is the key ingredient in the proof of Theorem \ref{t:logD} given in Section \ref{s:pfmain}. In Section \ref{s:Agroup} we determine the homotopy type of $Y_{A,k}$ for subgroups $A \leq G$ and comment on the optimality of the $\log D(G)$ factor in Theorem \ref{t:logD}. We conclude in Section \ref{s:con} with some remarks and open problems.

\section{Laplacians and their Eigenvalues}
\label{s:hodge}

Let $X$ be a finite simplicial complex on the vertex
set $V$. Let $X(k)$ denote the set of $k$-dimensional simplices in
$X$, each taken with an arbitrary but fixed orientation.
%The face numbers of $X$ are $f_k(X)=|X(k)|$.
A simplicial $k$-cochain is a real valued skew-symmetric function on
all ordered $k$-simplices of $X$. For $k \geq 0$ let $C^k(X)$
denote the space of $k$-cochains on $X$. The $i$-face of an
ordered $(k+1)$-simplex $\sigma=[v_0,\ldots,v_{k+1}]$ is the
ordered $k$-simplex
$\sigma_i=[v_0,\ldots,\widehat{v_i},\ldots,v_{k+1}]$. The
coboundary operator $d_k:C^k(X) \rightarrow C^{k+1}(X)$ is given
by $$d_k \phi (\sigma)=\sum_{i=0}^{k+1} (-1)^i \phi
(\sigma_i)~~.$$ It will be convenient to augment the cochain
complex $\{C^i(X)\}_{i \geq 0}$ with the $(-1)$-degree term
$C^{-1}(X)=\Com$ with the coboundary map $d_{-1}:C^{-1}(X)
\rightarrow C^0(X)$ given by $d_{-1}(a)(v)=a$ for $a \in \Com~,~v
\in V$. Let $Z^k(X)= \ker d_k$ denote the space of $k$-cocycles
and let $B^k(X)=\img d_{k-1}$ denote the space of
$k$-coboundaries. For $k \geq 0$ let $\thh^k(X)=Z^k(X)/B^k(X)~$
denote the $k$-th reduced cohomology group of $X$ with real
coefficients. For each $k \geq -1$ endow $C^k(X)$ with the
standard inner product $(\phi,\psi)_X=\sum_{\sigma \in X(k)}
\phi(\sigma)\psi(\sigma)~~$ and the corresponding $L^2$ norm
$||\phi||_X=(\phi,\phi)_X^{1/2}$. Let $d_k^*:C^{k+1}(X) \rightarrow C^k(X)$
denote the adjoint of $d_k$ with respect to these standard
inner products. The reduced $k$-th lower and upper Laplacians of $X$ are the positive semidefinite self-adjoint
maps of $C^k(X)$ given respectively by $L_k^-(X)=d_{k-1}d_{k-1}^*$ and $L_k^+(X)=d_k^*d_k$.
The $k$-th Laplacian of $X$ is $L_k(X)=L_k^-(X)+L_k^+(X)$.  Let $\mathcal{H}^k(X)=\ker L_k(X)=\ker d_{k-1}^* \cap \ker d_k$ denote the space of harmonic $k$-cochains. When there is no ambiguity concerning $X$, we shall abbreviate $L_k(X)=L_k$ and
$L_k^{\pm}(X)=L_k^{\pm}$. Clearly
\begin{equation*}
L_k^-\left(\img d_{k-1}\right) \subset \img d_{k-1}~~,~~
L_k^+\left(\img d_{k}^*\right) \subset \img d_{k}^*
\end{equation*}
\noindent
For a self-adjoint map $T$ on an inner product space $W$ let $\ccs(W,T)$ denote the
set of  eigenvalues of $T$ and let
$s(W,T,\lambda)$ denote the multiplicity of an eigenvalue $\lambda \in \ccs(W,T)$.
Let $\tilde{\ccs}(W,T)$ denote the multiset consisting of $s(W,T,\lambda)$ copies of each
$\lambda \in \ccs(W,T)$. The \emph{$k$-th spectral gap} of $X$ is
\[
\mu_k(X)=\min \ccs\left(C^k(X),L_k\right).
\]
\begin{remark}
\label{r:sgap}
(i) If $X$ is a graph then
$\mu_0(X)=\lambda_2(X)$. (ii) $\mu_k(X)>0$ iff $\thh^k(X;\Rea)=0$, hence $\mu_k$ may be viewed as a robustness measure of the property of having vanishing $k$-dimensional real cohomology.
\end{remark}
\noindent
The \emph{lower and upper $k$-th spectral gaps} of $X$ are defined respectively by
\[
\mu_k^-(X)=\min \ccs\left(\img d_{k-1},L_k^-\right)
\]
and
\[
\mu_k^+(X)=\min \ccs\left(\img d_{k}^*,L_k^+\right).
\]
\noindent
In Section \ref{s:sykg} we will use the some well known facts concerning Laplacians and their eigenvalues. For proofs, see e.g. \cite{DR02}.
\begin{proposition}
\label{p:hodge}
Let $0 \leq k \leq \dim X$. Then the following hold:
\begin{enumerate}[label=(\roman*)]

\item
Hodge Decomposition: There is an orthogonal direct sum decomposition:
\begin{equation}
\label{e:hodge1}
C^k(X)= \img d_{k-1} \oplus \mathcal{H}^k(X) \oplus \img d_{k}^*
\end{equation}

\item
\begin{equation}
\label{e:hodge3}
\ker L_k^-=\mathcal{H}^k(X) \oplus \img d_{k}^*~~~,~~~
\ker L_k^+=\img d_{k-1} \oplus \mathcal{H}^k(X).
\end{equation}

\item
Hodge isomorphism:
\begin{equation}
\label{e:hodge2}
\mathcal{H}^k(X) \cong \thh^k(X).
\end{equation}

\item For all $\lambda \neq 0$
\begin{equation}
\label{e:hodge4}
s\left(C^k(X),L_k,\lambda\right)=s\left(\img d_{k-1},L_k^-,\lambda\right)+
s\left(\img d_k^*,L_k^+,\lambda\right)
\end{equation}

\item
\begin{equation}
\label{e:hodge5}
\tilde{\ccs}\left(\img d_{k-1}^*,L_{k-1}^+\right)=\tilde{\ccs}\left(\img d_{k-1},L_{k}^-\right).
\end{equation}
\item If $H^k(X)=0$ then
\begin{equation*}
\label{e:hodge6}
\mu_k(X)=\min\left\{\mu_k^-(X),\mu_k^+(X)\right\}.
\end{equation*}
\end{enumerate}
\end{proposition}
\ \\ \\
In section \ref{s:sgyak} we shall use the following special case of Garland's fundamental eigenvalue estimate
(see Section 5 of \cite{Garland73} and Theorem 1.12 of \cite{BS97}).
The link of a simplex $\tau \in X(\ell)$ is
$X_{\tau}=\lk(X,\tau)=\{\eta \in X: \tau \cap \eta =\emptyset, \tau \cup \eta \in X\}$.
For $\phi \in C^j(X)$ and $\tau \in X(\ell)$ let $\phi_{\tau} \in C^{j-\ell-1}(X_{\tau})$ be defined by
$\phi_{\tau}(\eta)=\phi(\eta\tau)$, where $\eta\tau$ denotes the concatination of $\eta$ and $\tau$.
\begin{theorem}[Garland \cite{Garland73}]
\label{t:gar}
Let $X$ be a $k$-dimensional complex such that for all $\sigma \in X(k-1)$
\begin{equation*}
\deg_X(\sigma):=|\{\eta \in X(k): \sigma \subset \eta\}|=m.
\end{equation*}
Let $\lambda(X)=\min \left\{\lambda_2(X_{\tau}) : \tau \in X(k-2)\right\}$. Then
\begin{equation*}
\label{e:gar}
\min\left\{\frac{\|d_{k-1} \phi\|_X^2}{\|\phi\|_X^2}: 0 \neq \phi \in \ker d_{k-2}^*\right\}
\geq k\lambda(X)-(k-1)m.
\end{equation*}
\end{theorem}
\noindent
For completeness we indicate the proof. We first establish the following identity.
\begin{claim}
\label{c:garid}
For any $\phi \in C^{k-1}(X)$
\begin{equation}
\label{e:garid}
\|d_{k-1} \phi\|_X^2=\sum_{\tau \in X(k-2)} \|d_0 \phi_{\tau}\|_{X_{\tau}}^2-m(k-1) \|\phi\|_X^2.
\end{equation}
\end{claim}
\noindent
{\bf Proof.}
\begin{equation}
\label{e:gar1}
\begin{split}
&\|d_{k-1}\phi\|_X^2=\sum_{\sigma \in X(k)}|d_{k-1}\phi(\sigma)|^2 \\
&=\sum_{\sigma \in X(k)}\left(\sum_{i=0}^k (-1)^i \phi(\sigma_i)\right)
\left(\sum_{j=0}^k (-1)^j \phi(\sigma_j)\right) \\
&=\sum_{\sigma \in X(k)}\sum_{i=0}^k \phi(\sigma_i)^2 +2 \sum_{\sigma \in X(k)}\sum_{0 \leq i<j \leq k} (-1)^{i+j}
\phi(\sigma_i) \phi(\sigma_j) \\
&=m \|\phi\|_X^2-2 \sum_{\tau \in X(k-2)}\sum_{uv \in X_{\tau}(1)} \phi(u\tau)\phi(v\tau).
\end{split}
\end{equation}
\noindent
On the other hand
\begin{equation}
\label{e:gar2}
\begin{split}
&\sum_{\tau \in X(k-2)}\|d_0 \phi_{\tau}\|_{X_{\tau}}^2=
\sum_{\tau \in X(k-2)}
\sum_{uv \in X_{\tau}(1)}(\phi(v\tau)-\phi(u\tau))^2 \\
&=\sum_{\tau \in X(k-2)}\sum_{uv \in X_{\tau}(1)} \left(\phi(u\tau)^2+\phi(v\tau)^2\right) -
2 \sum_{\tau \in X(k-2)}\sum_{uv \in X_{\tau}(1)}
\phi(u\tau) \phi(v \tau) \\
&=m k \|\phi\|_X^2- 2 \sum_{\tau \in X(k-2)}\sum_{uv \in X_{\tau}(1)}
\phi(u\tau) \phi(v \tau).
\end{split}
\end{equation}
\noindent
Subtracting (\ref{e:gar2}) from (\ref{e:gar1}) we obtain (\ref{e:garid}).
{\enp}
\noindent
{\bf Proof of Theorem \ref{t:gar}.} Let $\phi \in \ker d_{k-2}^*$. Then for any $\tau \in X(k-2)$
\[
\sum_{v \in X_{\tau}(0)} \phi_{\tau}(v)=\sum_{v \in X_{\tau}(0)} \phi(v \tau)=d_{k-2}^*\phi(\tau)=0.
\]
Therefore, by the variational characterization of $\lambda_2(X_{\tau})$
\begin{equation}
\label{e:sgxt}
\|d_0\phi_{\tau}\|_{X_{\tau}}^2= (d_0^*d_0\phi_{\tau},\phi_{\tau})_{X_{\tau}} \geq \lambda_2(X_{\tau})\|\phi_{\tau}\|_{X_{\tau}}^2 \geq
\lambda(X)\|\phi_{\tau}\|_{X_{\tau}}^2.
\end{equation}
\noindent
Substituting (\ref{e:sgxt}) in (\ref{e:garid}) we obtain
\begin{equation*}
\begin{split}
\|d_{k-1} \phi\|_X^2&=\sum_{\tau \in X(k-2)} \|d_0 \phi_{\tau}\|_{X_{\tau}}^2-m(k-1) \|\phi\|_X^2 \\
&\geq \lambda(X) \sum_{\tau \in X(k-2)} \|\phi_{\tau}\|_{X_{\tau}}^2-m(k-1) \|\phi\|_X^2  \\
&= \big(\lambda(X) k- m(k-1)\big) \|\phi\|_X^2.
\end{split}
\end{equation*}
{\enp}

\section{Laplacians Spectra on $Y_{G,k}$}
\label{s:sykg}
In this section we prove the following characterization of the spectral gap of
complexes that contain the full $(k-1)$-skeleton of balanced complexes.
\begin{proposition}
\label{p:sgapy}
For any subcomplex $Y_{G,k}^{(k-1)} \subset Y \subset Y_{G,k}$
\begin{equation*}
\label{e:sgapy}
\mu_{k-1}(Y)
=\min\left\{\frac{\|d_{k-1}\phi\|_Y^2}{\|\phi\|_Y^2}:
0 \neq \phi \in \ker d_{k-2}^* \right\}.
\end{equation*}
\end{proposition}
We first record some facts concerning the Laplacian spectra of $Y_{G,k}$. We will use the notation introduced in Section \ref{s:hodge} with Laplacians $L_j=L_j(Y_{G,k})$.
\begin{proposition}
\label{p:ygk}
$~$
\begin{enumerate}[label=(\roman*)]
\item  For $0 \leq j \leq k$
\begin{equation}
\label{e:ygk1}
\begin{split}
\ccs\left(C^j(Y_{G,k}),L_j\right)&=\{tn: k-j \leq t \leq k+1\}, \\
s\left(C^j(Y_{G,k}),L_j,tn\right)&=\binom{k+1}{t} \binom{t}{k-j} \cdot (n-1)^{k+1-t}.
\end{split}
\end{equation}
\item  For $0 \leq j \leq k$
\begin{equation}
\label{e:ygk2}
\begin{split}
\ccs\left(\img d_{j-1},L_j^-\right)&=\{tn: k-j+1 \leq t \leq k+1\}, \\
s\left(\img d_{j-1},L_j^-,tn\right)&=\binom{k+1}{t} \binom{t-1}{k-j} \cdot (n-1)^{k+1-t}.
\end{split}
\end{equation}
\noindent
For $0 \leq j \leq k-1$
\begin{equation}
\label{e:ygk3}
\begin{split}
\ccs\left(\img d_j^*,L_j^+\right)&=\{tn: k-j \leq t \leq k+1\}, \\
s\left(\img d_j^*,L_j^+,tn\right)&=\binom{k+1}{t} \binom{t-1}{k-j-1} \cdot (n-1)^{k+1-t}.
\end{split}
\end{equation}
\end{enumerate}
\end{proposition}
\noindent
{\bf Proof.}
(i) Recall that $V_i$ is the $n$ point space $\{i\} \times G$.
For $0 \leq j \leq k$ let
\begin{equation*}
E_{k,j}=\left\{\ueps=
(\epsilon_1,\ldots,\epsilon_{k+1}) \in
\{-1,0\}^{k+1}: \epsilon_1+\cdots +\epsilon_{k+1}=j-k\right\}.
\end{equation*}
\noindent
The formula for the spectra of the Laplacians of joins
(see e.g. Section 4 in \cite{DR02}) implies that for $0 \leq j \leq k$
\begin{equation}
\label{e:splj}
\begin{split}
&\tilde{\ccs}\left(C^j(Y_{G,k}),L_j\right)= \tilde{\ccs}\left(C^j(V_1*\cdots*V_{k+1}),L_j\right) \\
&=\bigcup_{\ueps=(\epsilon_1,\ldots,\epsilon_{k+1}) \in E_{k,j}}
\left(\tilde{\ccs}\left(C^{\epsilon_1}(V_1),L_{\epsilon_1}\right)+\cdots+
\tilde{\ccs}\left(C^{\epsilon_{k+1}}(V_{k+1}),L_{\epsilon_{k+1}}\right)\right).
\end{split}
\end{equation}
\noindent
As $L_{-1}(V_i)$ is multiplication by $n$ and $L_0(V_i)$ is the all ones $n \times n$ matrix, it follows that
$\ccs(C^{-1}(V_i), L_{-1})=\{n\}$ and $\ccs(C^0(V_i),L_0)=\{0,n\}$ where
$s(C^0(V_i),L_0,0)=n-1$ and $s(C^0(V_i),L_0,n)=1$. Fix an
$\ueps=(\epsilon_1,\ldots,\epsilon_{k+1}) \in E_{k,j}$. Then
$I=\{1 \leq i \leq k+1: \epsilon_i=-1\}$ satisfies $|I|=k-j$.
The multiset corresponding to $\ueps$ in (\ref{e:splj}) is therefore
\begin{equation*}
\label{e:splj1}
M_{\ueps}=\overbrace{\{n\} +\cdots +\{n\}}^{k-j} +
\overbrace{\{0,\ldots,0,n\} +\cdots + \{0,\ldots,0,n\}}^{j+1}.
\end{equation*}
Clearly $M_{\ueps}$ consists of the elements $\{tn: k-j \leq t \leq k+1\}$,
where the multiplicity of $tn$ is $\binom{j+1}{t-(k-j)}(n-1)^{k+1-t}$.
Therefore
\begin{equation*}
\begin{split}
&s\left(C^j(Y_{G,k}),L_j,tn\right)=|E_{k,j}| \cdot \binom{j+1}{t-(k-j)}(n-1)^{k+1-t} \\
&=\binom{k+1}{k-j} \binom{j+1}{t-(k-j)}(n-1)^{k+1-t}
=\binom{k+1}{t} \binom{t}{k-j} \cdot (n-1)^{k+1-t}.
\end{split}
\end{equation*}
\\
(ii) We argue by decreasing induction on $j$. For the base case $j=k$, first note that (\ref{e:hodge3}) implies that
$0 \not\in \ccs(\img d_{k-1},L_k^-)$.
Moreover, as $L_k^+=0$ it follows by (\ref{e:hodge4}) that
for $\lambda \neq 0$
\begin{equation*}
s(\img d_{k-1},L_k^-,\lambda)=s(C^k(Y_{G,k}),L_k,\lambda).
\end{equation*}
Thus, by (\ref{e:ygk1})
\begin{equation*}
\ccs\left(\img d_{k-1},L_k^-\right)=\ccs\left(C^k(Y_{G,k}),L_k\right) \setminus \{0\}=
\{tn:1 \leq t \leq k+1\}
\end{equation*}
and
\begin{equation*}
s(\img d_{k-1},L_k^-,tn)=s(C^k(Y_{G,k}),L_k,tn)=\binom{k+1}{t}\cdot (n-1)^{k+1-t}.
\end{equation*}
For the induction step, let $1 \leq j_0 \leq k-1$ and assume that (\ref{e:ygk2}) holds for all
$j_0<j' \leq k$ and that (\ref{e:ygk3}) holds for all $j_0<j' \leq k-1$.
Then by (\ref{e:hodge5})
\begin{equation*}
\begin{split}
&\ccs\left(\img d_{j_0}^*,L_{j_0}^+\right)=\ccs\left(\img d_{j_0},L_{j_0+1}^-\right) \\
&=\{tn: k-(j_0+1)+1 \leq t \leq k+1\}=\{tn: k-j_0 \leq t \leq k+1\}
\end{split}
\end{equation*}
and
\begin{equation*}
\begin{split}
&s\left(\img d_{j_0}^*,L_{j_0}^+,tn\right)=s\left(\img d_{j_0},L_{j_0+1}^-,tn\right) \\
&= \binom{k+1}{t} \binom{t-1}{k-(j_0+1)} \cdot (n-1)^{k+1-t}=
\binom{k+1}{t} \binom{t-1}{k-j_0-1} \cdot (n-1)^{k+1-t}.
\end{split}
\end{equation*}
Thus (\ref{e:ygk3}) holds for $j=j_0$. Furthermore, by (\ref{e:hodge4})
\begin{equation*}
\begin{split}
&\{tn:k-j_0 \leq t \leq k+1\}= \ccs\left(C^{j_0}(X),L_{j_0}\right) \\
&=\ccs\left(\img d_{j_0-1},L_{j_0}^-\right)\cup
\ccs\left(\img d_{j_0}^*,L_{j_0}^+\right)
\end{split}
\end{equation*}
and for all $k-j_0 \leq t \leq k+1$
\begin{equation*}
\begin{split}
&s\left(\img d_{j_0-1},L_{j_0}^-,tn\right)=s\left(C^{j_0}(Y_{G,k}),L_{j_0},tn\right)
-s\left(\img d_{j_0}^*,L_{j_0}^+,tn\right) \\
&=\binom{k+1}{t} \binom{t}{k-j_0} \cdot (n-1)^{k+1-t}-
\binom{k+1}{t} \binom{t-1}{k-j_0-1} \cdot (n-1)^{k+1-t} \\
&=\binom{k+1}{t} \binom{t-1}{k-j_0} \cdot (n-1)^{k+1-t} \\
&=\left\{
\begin{array}{cl}
\binom{k+1}{t} \binom{t-1}{k-j_0} \cdot (n-1)^{k+1-t} & k-j_0+1 \leq t \leq k+1, \\
0 & t=k-j_0.
\end{array}
\right.
\end{split}
\end{equation*}
Thus (\ref{e:ygk2}) holds for $j=j_0$, thereby completing the inductive proof of (ii).
{\enp}
\noindent
{\bf Proof of Proposition \ref{p:sgapy}.} Let $Y_{G,k}^{(k-1)} \subset Y \subset Y_{G,k}$.
First note that the cases $j=k-1$ of (\ref{e:ygk3}) and (\ref{e:ygk2}) imply respectively that
\begin{equation*}
\label{e:vchr1}
\begin{split}
\alpha_+:&=\min \left\{\frac{\|d_{k-1} \phi\|_Y^2}{\|\phi\|_{Y}^2}: 0 \neq \phi \in \ker d_{k-2}^* \right\} \\
&\leq \min \left\{\frac{\|d_{k-1} \phi\|_{Y_{G,k}}^2}{\|\phi\|_{Y_{G,k}}^2}: 0 \neq \phi \in \ker d_{k-2}^* \right\} \\
&= \min \left\{\frac{\|d_{k-1} \phi\|_{Y_{G,k}}^2}{\|\phi\|_{Y_{G,k}}^2}:
0 \neq \phi \in d_{k-1}^* \left( C^k({Y_{G,k}})\right) \right\} \\
&=\mu_{k-1}^+({Y_{G,k}}) =\min\{ tn: 1 \leq t \leq k+1\}=n.
\end{split}
\end{equation*}
\noindent
and
\begin{equation*}
\label{e:vchr2}
\begin{split}
\alpha_-:&=\min \left\{\frac{\|d_{k-2}^* \phi\|_Y^2}{\|\phi\|_Y^2}: 0 \neq \phi \in \img d_{k-2} \right\} \\
&= \min \left\{\frac{\|d_{k-2}^* \phi\|_{Y_{G,k}}^2}{\|\phi\|_{Y_{G,k}}^2}: 0 \neq \phi \in \img d_{k-2} \right\} \\
&=\mu_{k-1}^-({Y_{G,k}})
=\min\{ tn: 2 \leq t \leq k+1\}=2n.
\end{split}
\end{equation*}
\noindent
Therefore $\alpha_+<\alpha_-$.
Moreover, $H^{k-1}({Y_{G,k}})=0$ together with (\ref{e:hodge1}) and (\ref{e:hodge2}) imply that there
is an orthogonal decomposition
\begin{equation*}
C^{k-1}(Y)=C^{k-1}({Y_{G,k}})=\img d_{k-2}\oplus \ker d_{k-2}^*.
\end{equation*}
Let $P_1,P_2$ denote the orthogonal projections of $C^{k-1}(Y)$ onto
$\img d_{k-2}$ and $\ker d_{k-2}^*$ respectively.
Then
\begin{equation*}
\label{e:vchr3}
\begin{split}
\mu_{k-1}(Y)&= \min\left\{\frac{(L_{k-1}\phi,\phi)_Y}{(\phi,\phi)_Y}:
0 \neq \phi \in C^{k-1}(Y)\right\} \\
&=\min\left\{\frac{\|d_{k-2}^*\phi\|_Y^2+\|d_{k-1}\phi\|_Y^2}{\|\phi\|_Y^2}:
0 \neq \phi \in C^{k-1}(Y)\right\} \\
&=\min\left\{\frac{\|d_{k-2}^*P_1\phi\|_Y^2+\|d_{k-1}P_2\phi\|_Y^2}{\|P_1 \phi\|_Y^2+\|P_2 \phi\|_Y^2 }:
0 \neq \phi \in C^{k-1}(Y)\right\} \\
&=\min \left\{\frac{\|d_{k-2}^*\phi_1\|_Y^2+\|d_{k-1}\phi_2\|_Y^2}{\|\phi_1\|_Y^2+\|\phi_2\|_Y^2}:
(0,0) \neq (\phi_1,\phi_2) \in \img d_{k-2} \times \ker d_{k-2}^*\right\} \\
&=\min\{\alpha_-,\alpha_+\}=\alpha_+.
\end{split}
\end{equation*}
{\enp}

\section{The Fourier Transform}
\label{s:fourier}

Let $\cll(G)$ denote the algebra of complex valued functions on $G$ with the convolution product
$\phi*\psi (x)=\sum_{y \in G} \phi(y)\psi(y^{-1}x)$.
The inner product on $\cll(G)$ is given by
\[
\langle\phi,\psi\rangle=\sum_{x \in G} \phi(x) \overline{\psi(x)}.
\]
The Frobenius inner product and norm on $M_d(\Com)$ are given respectively by
$\langle S,T \rangle=\tr(ST^*)$ and $\|T\|_F=\sqrt{\langle T,T \rangle}=\sqrt{\tr(TT^*)}$. The Frobenius norm of a product satisfies
\begin{equation}
\label{e:submul}
\|ST\|_F \leq \|S\|\cdot \|T\|_F.
\end{equation}
\noindent
Let $\crr(G)$ denote the algebra $\prod_{\rho \in \ghat} M_{d_{\rho}}(\Com)$ with coordinate wise addition and multiplication. Define an inner product on $\crr(G)$ by
\[
\left\langle \left(S_{\rho}: \rho \in \ghat\right),\left(T_{\rho}: \rho \in \ghat\right)\right\rangle=
\frac{1}{n}\sum_{\rho} d_{\rho}\langle S_{\rho}, T_{\rho}\rangle=
\frac{1}{n}\sum_{\rho} d_{\rho}\tr(S_{\rho} T_{\rho}^*).
\]
The associated norm is given by
\[
\left\|\left(T_{\rho}: \rho \in \ghat\right)\right\|_F=\left(\frac{1}{n} \sum_{\rho \in \ghat} d_{\rho} \|T_{\rho}\|_F^2 \right)^{\frac{1}{2}}.
\]

\begin{definition}
\label{d:ft}
For $\phi \in \cll(G)$ and a representation $\rho$ of $G$ of degree $d$ let
\[
\phat(\rho)=\sum_{x \in G} \phi(x) \rho(x) \in M_{d}(\Com).
\]
\noindent
The {\bf Fourier Transform} $\cf:\cll(G) \rightarrow \crr(G)$ is given by
\[
\cf(\phi)=\left(\phat(\rho): \rho \in \ghat \right).
\]
\end{definition}
\noindent
A basic result in representation theory (see e.g. exercise 3.32 in \cite{Fulton-Harris91}) asserts that $\cf$ is an isomorphism of algebras and an isometry. In particular, $\cf$ satisfies the Parseval identity: For any $\phi,\psi \in \cll(G)$
\begin{equation}
\label{e:parseval}
\langle \phi,\psi\rangle=\left\langle \cf(\phi),\cf(\psi)\right\rangle
=
\frac{1}{n} \sum_{\rho \in \ghat} d_{\rho} \left\langle \widehat{\phi}(\rho),
\widehat{\psi}(\rho)\right\rangle.
\end{equation}

\section{The $(k-1)$-Spectral Gap of $Y_{A,k}$}
\label{s:sgyak}

In this section we prove Theorem \ref{t:specgap}. Let $X=Y_{A,k}$. We need two preliminary observations.
Let $\rmc_A$ be the graph on the vertex set $V(\rmc_A)=\{1,2\} \times G$ with edge set
\[
E(\rmc_A)=\left\{\{(1,x_1),(2,x_2)\}: x_1\cdot x_2 \in A\right\}.
\]
\begin{claim}
\label{c:linktau}
For any $\tau \in X(k-2)$, the graph $X_{\tau}=\lk(X,\tau)$ is isomorphic to
$\rmc_A$.
\end{claim}
\noindent
{\bf Proof.} Let $\tau=\{(j,y_j)\}_{j \in J}$ where $J \subset [k+1]:=\{1,\ldots,k+1\}$ and $|J|=k-1$.
Let $[k+1]\setminus J=\{i_1<i_2\}$.
Let $z_1=y_1 \cdots y_{i_1-1}$, $z_2=y_{i_1+1} \cdots y_{i_2-1}$ and  $z_3=y_{i_2+1} \cdots y_{k+1}$.
Then $X_{\tau}$ is the graph on the vertex set
$V_{\tau}=\{i_1,i_2\} \times G$ with edge set
\[
E_{\tau}=\left\{ \{(i_1,x_{i_1}),(i_2,x_{i_2})\}: z_1 x_{i_1} z_2 x_{i_2} z_3 \in A\right\}.
\]
Let $\varphi: V_{\tau} \rightarrow V(\rmc_A)$ be given by
\[
\varphi\left((i_t,x_{i_t})\right)=
\left\{
\begin{array}{ll}
(1,z_1x_{i_1}z_2) & t=1, \\
(2,x_{i_2}z_3) & t=2.
\end{array}
\right.
\]
Then $\varphi$ is an isomorphism between $X_{\tau}$ and $\rmc_A$.
{\enp}
\noindent
The next result gives a lower bound on the spectral gap of $\rmc_A$.
\begin{proposition}
\label{p:sggr}
\begin{equation*}
\label{e:sggr}
\lambda_2(\rmc_A) \geq |A|-\nu(A).
\end{equation*}
\end{proposition}
\noindent
{\bf Proof.} Let $\phi \in C^0\left(V(\rmc_A)\right)$ such that $\sum_{v \in V(\rmc_A)} \phi(v)=0$.
For $i=1,2$ let $\phi_i \in \cll(G)$ be given by $\phi_i(x)=\phi\left((i,x)\right)$.
Define $\psi \in \cll(G)$
by $\psi(x)=\phi_2(x^{-1})$ and for $a \in A$ let $\psi_a(x)=\psi(a^{-1} x)=\phi_2(x^{-1}a)$.
Then
\begin{equation*}
\begin{split}
&\widehat{\phi_1}(\one)+\widehat{\psi_a}(\one)=\left(\sum_{x \in G} \phi_1(x)\right)
+\left(\sum_{x \in G} \psi_a(x)\right) \\
&=\left(\sum_{x \in G} \phi_1(x)\right)+\left(\sum_{x \in G} \phi_2(x)\right)=
\sum_{v \in V(\rmc_A)} \phi(v)=0.
\end{split}
\end{equation*}
Hence
\begin{equation}
\label{e:psia0}
\widehat{\phi_1}(\one) \cdot \widehat{\psi_a(\one})=-\widehat{\phi_1}(\one)^2 \leq 0.
\end{equation}
\noindent
For any $\rho \in \ghat$
\begin{equation}
\label{e:psia1}
\begin{split}
\widehat{\psi_a}(\rho)&=\sum_{x \in G} \psi_a(x)\rho(x)
=\sum_{x \in G} \phi_2(x^{-1}a) \rho(x) \\
&= \sum_{y \in G} \phi_2(y) \rho(a y^{-1})
=\rho(a) \sum_{y \in G} \phi_2(y) \rho(y^{-1})=\rho(a) \widehat{\psi}(\rho).
\end{split}
\end{equation}
\noindent
Using the Parseval identity (\ref{e:parseval}) together with (\ref{e:psia0}), (\ref{e:psia1}) and (\ref{e:submul}) we obtain
\begin{equation}
\label{e:psia2}
\begin{split}
&\sum_{a \in A} \left\langle \phi_1,\psi_a\right\rangle=\sum_{a \in A} \left\langle \cf(\phi_1),\cf(\psi_a)\right\rangle \\
&=\frac{1}{n} \sum_{a \in A} \left(
\sum_{\rho \in  \ghat} d_{\rho} \left\langle \widehat{\phi_1}(\rho), \widehat{\psi_a}(\rho)\right\rangle\right)
\\
&=\frac{1}{n}\sum_{a \in A} \left( \widehat{\phi_1}(\one) \cdot \widehat{\psi_a}(\one) +
\sum_{\rho \in  \ghat_+} d_{\rho} \left\langle \widehat{\phi_1}(\rho), \widehat{\psi_a}(\rho)\right\rangle\right)
\\
&=-\frac{|A|}{n}\widehat{\phi_1}(\one)^2+
\frac{1}{n} \sum_{a \in A}\sum_{\rho \in \ghat_+} d_{\rho}\left\langle\widehat{\phi_1}(\rho), \rho(a)\widehat{\psi}(\rho)\right\rangle
\\
&=-\frac{|A|}{n}\widehat{\phi_1}(\one)^2+\frac{1}{n}\sum_{\rho \in \ghat_+}d_{\rho}\left\langle\widehat{\phi_1}(\rho), \widehat{1_A}(\rho) \widehat{\psi}(\rho)\right\rangle \\
&\leq \frac{1}{n}\sum_{\rho \in \ghat_+}d_{\rho}\|\widehat{\phi_1}(\rho)\|_F \cdot \|\widehat{1_A}(\rho)\cdot  \widehat{\psi}(\rho)\|_F \\
&\leq \frac{1}{n}\sum_{\rho \in \ghat_+}d_{\rho}\|\widehat{\phi_1}(\rho)\|_F \cdot \|\widehat{1_A}(\rho)\| \cdot  \|\widehat{\psi}(\rho)\|_F \\
&\leq \nu(A) \sum_{\rho\in \ghat} \left( \sqrt{\frac{d_{\rho}}{n}} \|\widehat{\phi_1}(\rho)\|_F \right) \cdot
\left( \sqrt{\frac{d_{\rho}}{n}} \|\widehat{\psi}(\rho)\|_F \right) \\
&\leq \nu(A) \left(\frac{1}{n}\sum_{\rho\in \ghat} d_{\rho} \|\widehat{\phi_1}(\rho)\|_F^2 \right)^{\frac{1}{2}}
\left(\frac{1}{n}\sum_{\rho\in \ghat} d_{\rho} \|\widehat{\psi}(\rho)\|_F^2 \right)^{\frac{1}{2}} \\
&=\nu(A)\cdot \|\phi_1\|\cdot \|\psi\| =\nu(A)\cdot \|\phi_1\|\cdot \|\phi_2\|.
\end{split}
\end{equation}
\noindent
Finally by (\ref{e:psia2})
\begin{equation*}
\label{e:dephi}
\begin{split}
\|d_0 \phi\|_{\rmc_A}^2&=\sum_{\{(x_1,x_2) \in G^2: x_1x_2 \in A\}} d_0\phi \left( [(1,x_1),(2,x_2)]\right)^2 \\
&=\sum_{\{(x_1,x_2) \in G^2: x_1x_2 \in A\}} \left(\phi_2(x_2)-\phi_1(x_1)\right)^2 \\
&=\sum_{a \in A} \sum_{x \in G} \left(\phi_2(x^{-1}a)-\phi_1(x)\right)^2 \\
&=\sum_{a \in A} \sum_{x \in G} \left( \phi_1(x)^2+\phi_2(x^{-1}a)^2-2\phi_1(x)\psi_a(x)\right)\\
&=|A| \left(\|\phi_1\|^2+\|\phi_2\|^2\right)-2\sum_{a \in A} \langle \phi_1,\psi_a\rangle \\
&\geq |A| \cdot \|\phi\|_{\rmc_A}^2-2 \nu(A)\cdot \|\phi_1\|\cdot \|\phi_2\| \\
&\geq |A| \cdot \|\phi\|_{\rmc_A}^2 -\nu(A) \left(\|\phi_1\|^2 + \|\phi_2\|^2 \right) \\
&= \left(|A|-\nu(A)\right) \|\phi\|_{\rmc_A}^2.
\end{split}
\end{equation*}
Therefore $\lambda_2(\rmc_A) \geq |A|-\nu(A)$.
{\enp}
\noindent
{\bf Proof of Theorem \ref{t:specgap}.} Clearly $\deg_X(\sigma)=|A|$ for all $\sigma \in X(k-1)$.
By Claim \ref{c:linktau} and Proposition \ref{p:sggr}
\begin{equation}
\label{e:blamx}
\lambda(X)=\min\left\{\lambda_2(X_{\tau}): \tau \in X(k-2)\right\} = \lambda_2(\rmc_A) \geq |A|-\nu(A).
\end{equation}
Using Proposition \ref{p:sgapy}, Garland's Theorem \ref{t:gar} and (\ref{e:blamx}) it follows that
\begin{equation*}
\begin{split}
\mu_{k-1}(X)&=\min\left\{\frac{\|d_{k-1} \phi\|_X^2}{\|\phi\|_X^2}: 0 \neq \phi \in \ker d_{k-2}^*\right\} \\
&\geq k\lambda(X)-(k-1)|A| \geq k (|A|-\nu(A))-(k-1)|A| \\
&=|A|-k\cdot\nu(A).
\end{split}
\end{equation*}
{\enp}

\section{The Spectral Gap of a Random $Y_{A,k}$}
\label{s:pfmain}

In this section we prove Theorem \ref{t:logD}.
We will use the following matrix version of Bernstein's large deviation inequality due to
Tropp (Theorem 1.6 in \cite{Tropp12}).
\begin{theorem}[\cite{Tropp12}]
\label{t:tropp}
Let $\{X_i\}_{i=1}^m$ be independent random variables taking values in $M_d(\Com)$ such that
$E[X_i]=0$ and $\|X_i\| \leq R$ for all $1 \leq i \leq m$, and let
\[
\sigma^2=\max\left\{\left\|\sum_{i=1}^m E[X_iX_i^*]\right\|,\left\|\sum_{i=1}^m E[X_i^*X_i]\right\|\right\}.
\]
Then for any $\lambda\geq 0$
\begin{equation*}
\label{e:mjcft}
\prob\left[\left\|\sum_{i=1}^m X_i\right\| \geq \lambda\right] \leq 2d\exp\left(-\frac{3\lambda^2}{6\sigma^2+2R\lambda}\right).
\end{equation*}
\end{theorem}
\noindent
{\bf Proof of Theorem \ref{t:logD}.}
Let $k \geq 1$ and $0<\epsilon<1$ be fixed and let $m=\lceil 9 k^2\log D(G)/\epsilon^2 \rceil$. Let
$\Omega$ denote the uniform probability space of all $m$-subsets of $G$.
Suppose that $A \in \Omega$ satisfies $\nu(A) \leq \epsilon k^{-1}m$. Then by Theorem \ref{t:specgap}
\begin{equation*}
\begin{split}
\mu_{k-1}(X_{A,k}) &\geq |A|-k \cdot \nu(A) \\
& \geq m-k \cdot \epsilon k^{-1} m=(1-\epsilon)m.
\end{split}
\end{equation*}
Theorem \ref{t:logD} will therefore follow from
\begin{proposition}
\label{p:cher}
\begin{equation*}
\label{e:cher}
\pr_{\Omega}\left[~A \in \Omega~:~\max_{\rho \in \ghat_+}
\|\widehat{1_A}(\rho)\|>\epsilon k^{-1} m~\right]< \frac{6}{n}.
\end{equation*}
\end{proposition}
\noindent
{\bf Proof.}
Let $\rho \in \ghat_+$ be fixed and let $\lambda=\epsilon k^{-1} m$.
Let $\Omega'$ denote the uniform probability space $G^m$, and for $1 \leq i  \leq m$ let $X_i$ be the random variable defined on $\omega'=(a_1,\ldots,a_m) \in \Omega'$ by
$X_i(\omega')=\rho(a_i) \in U(d_{\rho})$. As $\rho \in \ghat_+$, it follows by Schur's Lemma that $E[X_i]=0$.
It is straightforward to check the $X_i$'s also satisfy the additional conditions of Theorem \ref{t:tropp} with $\sigma^2=m$ and $R=1$. Hence
\begin{equation}
\label{e:bernstein}
\begin{split}
&\pr_{\Omega'}\left[~\omega' \in \Omega'~:~\left\| \sum_{i=1}^m X_i(\omega') \right\| \geq \lambda~\right]
\leq 2d_{\rho}\exp\left(-\frac{3\lambda^2}{6\sigma^2+2R\lambda}\right) \\
&=2d_{\rho} \exp\left(-\frac{3(\epsilon k^{-1}m)^2}{6m+2\epsilon k^{-1}m}\right)
< 2d_{\rho} \exp\left(-\frac{\epsilon^2 m}{3k^2}\right) \\
&\leq 2d_{\rho} \exp\left(-3 \log D(G)\right)=2d_{\rho}D(G)^{-3}.
\end{split}
\end{equation}
Let $\Omega''=\{(a_1,\ldots,a_m) \in G^m: a_i \neq a_j \text{~for~} i \neq j \}$
denote the subspace of $\Omega'$ consisting of all sequences in $G^m$ with pairwise distinct elements. Note that the assumption $n > 10^6 k^8 \epsilon^{-8}$ implies that
$\frac{m^2}{n-m} <1$
and therefore
\begin{equation}
\label{e:omin}
\pr_{\Omega'}[~\Omega''~]=\prod_{i=1}^m \frac{n-i+1}{n}
> \left(\frac{n-m}{n}\right)^m \geq \exp\left(-\frac{m^2}{n-m}\right) \geq e^{-1}.
\end{equation}
\noindent
Combining (\ref{e:bernstein}) and (\ref{e:omin}) we obtain
\begin{equation}
\label{e:bernstein1}
\begin{split}
&\pr_{\Omega}\left[A \in \Omega:\|\widehat{1_A}(\rho)\| \geq \epsilon k^{-1} m~\right] \\
&=\pr_{\Omega''}\left[~\omega'' \in \Omega'':\left\|\sum_{i=1}^m X_i(\omega'')\right\|
\geq \epsilon k^{-1} m~\right] \\
&\leq \pr_{\Omega'}\left[~\omega' \in \Omega'~:~\left\|\sum_{i=1}^m X_i(\omega')\right\| \geq
\epsilon k^{-1} m~\right]
\cdot \left(\pr_{\Omega'}[~\Omega''~]\right)^{-1} \\
&< 6d_{\rho}D(G)^{-3}.
\end{split}
\end{equation}
Note that
\begin{equation*}
D(G)^2=\left(\sum_{\rho \in \ghat} d_{\rho}\right)^2 \geq  \sum_{\rho \in \ghat} d_{\rho}^2 =n.
\end{equation*}
\noindent
Using the union bound and (\ref{e:bernstein1}) it thus follows that
\begin{equation*}
\pr_{\Omega}\left[~\nu(A)\geq \epsilon k^{-1}m~\right]
<6\sum_{\rho \in \ghat} d_{\rho} D(G)^{-3} = 6 D(G)^{-2}< \frac{6}{n}.
\end{equation*}
{\enp}

\section{$Y_{A,k}$ for a Subgroups $A$}
\label{s:Agroup}

Let $A$ be a subgroup of $G$ of order $|A|=m$ and let $\ell=(G:A)=\frac{n}{m}$. Let
\[
\gamma_0(m,k)=
(n-m)n^k+\ell^k(m-1)^{k+1}-(n-1)^{k+1}
\]
and
\[
\gamma_1(m,k)=\ell^k(m-1)^{k+1}.
\]
\noindent
The homotopy type of $Y_{A,k}$ is given by the following
\begin{proposition}
\label{p:hoty}
(i)
\begin{equation}
\label{e:hoty1}
Y_{A,1}\simeq \displaystyle{\coprod^{\ell} \bigvee^{(m-1)^2}\sph^1}.
\end{equation}
\noindent
(ii) For $k \geq 2$
\begin{equation}
\label{e:hoty2}
Y_{A,k}\simeq
\displaystyle{\bigvee^{\gamma_0(m,k)} \sph^{k-1} \vee \bigvee^{\gamma_1(m,k)} \sph^k}.
\end{equation}
\end{proposition}
\noindent
The proof of Proposition \ref{p:hoty}(ii) depends on the Wedge Lemma of
Ziegler and \v{Z}ivaljevi\'{c} (Lemma 1.8 in \cite{ZZ93}). The version below appears in
\cite{HRW98}. For a poset $(P,\prec)$ and $p \in P$ let $P_{\prec p}=\{q \in P: q \prec p\}$. Let $\Delta(P)$ denote the order complex of $P$. Let $Y$ be a regular CW-complex and let $\{Z_i\}_{i=1}^{\ell}$ be subcomplexes of $Y$ such that $\bigcup_{i=1}^{\ell} Z_i=Y$. Let $(P,\prec)$ be the poset whose elements index all distinct partial intersections $\bigcap_{j \in J} Z_j$, where $\emptyset \neq J \subset [\ell]$. Let $U_p$ denote the partial intersection indexed by $p \in P$,
and let $\prec$ denote reverse inclusion, i.e. $p \prec q$ if $U_q \subsetneqq U_p$.

\vspace{2mm}
\noindent
{\bf Wedge Lemma \cite{ZZ93,HRW98}.}
suppose that for any $p \in P$ there exists a $c_p \in U_p$ such that the inclusion
$\bigcup_{q \succ p} U_q \hookrightarrow U_p$ is homotopic to the constant map to $c_p$. Then
\begin{equation*}
\label{e:wlemma}
Y \simeq \bigvee_{p \in P} \Delta(P_{\prec p})*U_p.
\end{equation*}
\noindent
{\bf Proof of Proposition \ref{p:hoty}.} Let
$g_1,\ldots,g_{\ell} \in G$ be coset representatives of $A$, i.e. $G=\bigcup_{i=1}^{\ell} g_i A$.
(i) The graph $Y_{A,1}$ is isomorphic to the disjoint union
$\coprod_{i=1}^{\ell} Ag_i^{-1} * g_iA$. This implies (\ref{e:hoty1}) since each $Ag_i^{-1} * g_iA$ is a complete $m$ by $m$ bipartite graph
and hence homotopic to  a wedge of
$(m-1)^2$ circles.
\\
(ii) Let $k\geq 2$. For $1 \leq i \leq \ell$ let
$W_{k,i}=\{k+1\} \times g_i A  \subset V_{k+1}$ and let
\begin{equation}
\label{e:zizi}
Z_{k,i}=Y_{Ag_i^{-1},k-1}*W_{k,i} \cong Y_{A,k-1}*[m].
\end{equation}
\noindent
Then $\bigcup_{i=1}^{\ell} Z_{k,i}=Y_{A,k}$. Indeed,
let $x_1,\ldots,x_{k+1} \in G$, and suppose that $x_{k+1} \in g_i A$. Then
\begin{equation*}
\begin{split}
\sigma=\{(1,x_1),\ldots,(k,x_k),(k+1,x_{k+1})\} \in Y_{A,k} &~\Longleftrightarrow~  x_1\cdots x_{k+1} \in A \\
~\Longleftrightarrow~ x_1\cdots x_k \in Ag_i^{-1} &~\Longleftrightarrow~ \sigma \in Z_{k,i}(k).
\end{split}
\end{equation*}
Moreover, for any $1 \leq j \neq j' \leq t$
\begin{equation}
\label{e:intzij}
Z_{k,j} \cap Z_{k,j'} =\bigcap_{i=1}^{\ell} Z_{k,i}=Y_{G,k-1}^{(k-2)}.
\end{equation}
Let $N_k=(-1)^{k-2}\tilde{\chi}\left(Y_{G,k-1}^{(k-2)}\right)=n^k-(n-1)^k$.
As $Y_{G,k-1}^{(k-2)}$ is a matroidal complex of rank $k-1$, it follows (see e.g. Theorem 7.8.1 in \cite{Bjorner92}) that
\begin{equation}
\label{e:valofn}
Y_{G,k-1}^{(k-2)} \simeq
\bigvee^{N_k} \sph^{k-2}.
\end{equation}
Eq. (\ref{e:intzij}) implies that the intersection poset $(P,\prec)$ of the cover $\{Z_{k,i}\}_{i=1}^{\ell}$ is $P=[\ell] \cup \{\widehat{1}\}$, where $i \in [\ell]$ represents $Z_{k,i}$, $\widehat{1}$ represents $Y_{G,k-1}^{(k-2)}$,
$[\ell]$ is an antichain and
$i \prec \widehat{1}$ for all $i \in [\ell]$.
Note that
$\Delta(P_{\prec i})=\emptyset$ for all $i \in [m]$ and $\Delta(P_{\prec \widehat{1}})$ is the discrete space $[\ell]$. We proceed to prove (\ref{e:hoty2}) by induction on $k$. We first establish the induction step.
Let $k \geq 3$ and assume that (\ref{e:hoty2}) holds for $k-1$.
Then $Z_{k,i} \cong Y_{A,k-1}*[m]$
 is homotopy equivalent to a wedge of spheres of dimensions $k-1$ and $k$. As $Y_{G,k-1}^{(k-2)}$ is a wedge of $(k-2)$-spheres, it follows that
the inclusion $Y_{G,k-1}^{(k-2)} \hookrightarrow Z_{k,i}$ is null homotopic.
Applying the Wedge Lemma together with (\ref{e:zizi}), (\ref{e:valofn}) and the induction hypothesis, we obtain
\begin{equation}
\label{e:wedge1}
\begin{split}
Y_{A,k} &\simeq
\left( \bigvee_{i \in [\ell]} \Delta(P_{\prec i}) * Z_{k,i}\right) \vee
\left( \Delta(P_{\prec \widehat{1}}) * Y_{G,k-1}^{(k-2)} \right) \\
&=\left(\bigvee_{i \in [\ell]}Z_{k,i}\right) \vee \left([\ell]*Y_{G,k-1}^{(k-2)}\right) \\
&\cong\left(\bigvee_{i \in [\ell]} Y_{A,k-1}* [m]\right) \vee  \left([\ell]*\bigvee^{N_k} \sph^{k-2}\right) \\
&\simeq
\left(\bigvee_{i \in [\ell]} \left(\displaystyle{\bigvee^{\gamma_0(m,k-1)} \sph^{k-2} \vee \bigvee^{\gamma_1(m,k-1)} \sph^{k-1}}\right) * [m]\right) \vee  \left([\ell]*\bigvee^{N_k} \sph^{k-2}\right) \\
&\simeq\bigvee_{i \in [\ell]} \left(\displaystyle{\bigvee^{(m-1)\gamma_0(m,k-1)} \sph^{k-1} \vee \bigvee^{(m-1)\gamma_1(m,k-1)} \sph^{k}}\right) \vee  \left(\bigvee^{(\ell-1)N_k} \sph^{k-1}\right) \\
&\simeq\bigvee^{t_0} \sph^{k-1} \vee \bigvee^{t_1} \sph^{k},
\end{split}
\end{equation}
where
\begin{equation*}
\begin{split}
t_0&=\ell(m-1)\gamma_0(m,k-1)+(\ell-1)N_k \\
&=\ell(m-1)\left((n-m)n^{k-1}+\ell^{k-1}(m-1)^{k}-(n-1)^{k}\right)+
(\ell-1)\left(n^k-(n-1)^k\right) \\
&=(n-m)n^k+\ell^k(m-1)^{k+1}-(n-1)^{k+1}=\gamma_0(m,k)
\end{split}
\end{equation*}
and
\begin{equation*}
\begin{split}
t_1&=\ell(m-1)\gamma_1(m,k-1)=\ell(m-1)\left(\ell^{k-1}(m-1)^{k}\right) \\
&=\ell^{k}(m-1)^{k+1}=\gamma_1(m,k).
\end{split}
\end{equation*}
\noindent
This completes the induction step. To prove (\ref{e:hoty2}) for $k=2$, first note that
assumptions of the Wedge lemma hold for the decomposition $Y_{A,2}= \bigcup_{i=1}^{\ell} Z_{2,i}$.
Arguing as in (\ref{e:wedge1}), it thus follows that
\begin{equation*}
\label{e:wedge2}
\begin{split}
Y_{A,2} &\simeq
\left(\bigvee_{i \in [\ell]} Y_{A,1}* [m]\right) \vee  \left([\ell]*\bigvee^{N_2} \sph^{0}\right) \\
&\simeq \left(\bigvee_{i \in [\ell]} \left(\displaystyle{\coprod^{\ell} \bigvee^{(m-1)^2}\sph^1}\right)* [m]\right) \vee  \left([\ell]*\bigvee^{2n-1} \sph^{0}\right) \\
&\simeq \bigvee^{\ell} \bigvee^{m-1} \left(\bigvee^{\ell(m-1)^2}\sph^2 \vee \bigvee^{\ell-1} \sph^1\right) \vee  \bigvee^{(\ell-1)(2n-1)} \sph^{1} \\
&=\bigvee^{t_0} \sph^{1} \vee \bigvee^{t_1} \sph^{2},
\end{split}
\end{equation*}
where
\begin{equation*}
\begin{split}
t_0&=\ell(m-1)(\ell-1)+(\ell-1)(2n-1) \\
&=(n-m)n^2+\ell^2(m-1)^{3}-(n-1)^{3}=\gamma_0(m,2)
\end{split}
\end{equation*}
and
\begin{equation*}
\begin{split}
t_1&=\ell^2(m-1)^3=\gamma_1(m,k).
\end{split}
\end{equation*}
\noindent
This completes the proof of the base case $k=2$ and of the Proposition.
{\enp}
\noindent
As $\gamma_0(m,k)>0$ for all $m<n$, it follows from Proposition \ref{p:hoty} that if $A \subset G$
generates a subgroup $\langle A \rangle$ of order $m<n$ then
\begin{equation*}
\tilde{\beta}_{k-1}(Y_{A,k}) \geq \tilde{\beta}_{k-1}(Y_{\langle A \rangle,k})
= \gamma_0(m,k)>0
\end{equation*}
and therefore $\mu_{k-1}(Y_{A,k})=0$.
This implies that the $\log D(G)=\Theta(\log n)$ factor in Theorem \ref{t:logD} cannot in general be improved.

\section{Concluding Remarks}
\label{s:con}

In this paper we studied the $(k-1)$-spectral gap of complexes $Y_{A,k}$ where $A$ is a subset of a finite group $G$. Our main results included a lower bound on $\mu_{k-1}(Y_{A,k})$ in terms of the Fourier transform of $1_A$
and a proof that for a sufficiently large constant $c(k,\epsilon)$, if $A$ is a random subset of $G$ of size at least $c(k,\epsilon) \log D(G)$, then $Y_{A,k}$ has a nearly optimal $(k-1)$-th spectral gap.
In view of Remark \ref{r:sgap}(ii) it
would be interesting to find suitable counterparts of Theorems \ref{t:specgap} and \ref{t:logD} for other robustness measures of cohomological triviality, e.g. for coboundary expansion.
We briefly recall the relevant definitions.
For a simplicial complex $X$ and a binary $k$-cochain $\phi \in C^{k}(X;\FF_2)$, let
\[
\|\phi\|_{\rm H}=\left|\left\{\sigma \in X(k): \phi(\sigma) \neq 0\right\}\right|
\]
denote the Hamming norm of  $\phi$ and
let
\[
\|\phi\|_{\csy}=\min \left\{|\supp(\phi+d_{k-1}\psi)|: \psi \in C^{k-1}(X;\FF_2)\right\}
\]
denote the cosystolic norm of $\phi$. The $k$-th coboundary expansion constant of $X$ (see e.g.
\cite{Lubotzky18}) is given by
\[
h_{k}(X)=\min\left\{\frac{\|d_k\phi\|_{{\rm H}}}{\|\phi\|_{\csy}}:
\phi \in C^k(X;\FF_2) \setminus B^k(X;\FF_2)\right\}.
\]
\noindent
In light of Theorem \ref{t:logD} we suggest the following
\begin{conjecture}
\label{c:logDc}
For any fixed $k \geq 1$ there exist constants $C(k)< \infty$ and $\epsilon(k)>0$ such that for any group $G$,
the random balanced Cayley complexes $Y_{A,k}$ with $|A|=C(k) \log D(G)$ satisfy $h_{k-1}(Y_{A,k}) \geq \epsilon(k)$ a.a.s. as $|G| \rightarrow \infty$.
\end{conjecture}

In a different direction, consider the following example of balanced Cayley complexes.
Let $p,q$ be distinct odd primes such that $q>2 \sqrt{p}$ and $\left(\frac{p}{q}\right)=1$, and let
$G_q=PSL_2(\FF_q)$. The celebrated construction of Ramanujan graphs
by Lubotzky, Phillips and Sarnak \cite{LPS88} implies that there exists a subset
$S_{p,q} \subset G_q$ of cardinality $|S_{p,q}|=p+1$ such that
$\nu(S_{p,q}) \leq 2 \sqrt{p}$.
If $p \geq 4k^2$ then by Theorem \ref{t:specgap}
\begin{equation}
\label{e:yspq}
\mu_{k-1}\left(Y_{S_{p,q},k}\right) \geq |S_{p,q}|-k\cdot \nu\left(S_{p,q}\right)
\geq (p+1)-2k\sqrt{p}\geq 1.
\end{equation}
\noindent
The following conjecture may be viewed as a coboundary expansion analogue of (\ref{e:yspq}).
\begin{conjecture}
\label{c:ramanujan}
For any fixed $k \geq 1$ there exist constants $p_0(k)<\infty$ and $\epsilon_0(k)>0$ such that
if $p>p_0(k)$, $q>2 \sqrt{p}$ and $\left(\frac{p}{q}\right)=1$ then
$h_{k-1}\left(Y_{S_{p,q},k}\right) \geq \epsilon_0(k)$.
\end{conjecture}


\begin{thebibliography}{99}

\bibitem{AR94}
N. Alon and Y. Roichman,
Random Cayley graphs and expanders,
{\it Random Structures Algorithms}, {\bf 5}(1994) 271-–284.

\bibitem{BS97}
W. Ballmann and J. \'{S}wi\c{a}tkowski, On $L^2$-cohomology and property
(T) for automorphism groups of polyhedral cell complexes, {\it
Geometric and Functional Analysis}, {\bf 7}(1997) 615-645.

\bibitem{BAM20}
O. Beit-Aharon and R. Meshulam, Spectral expansion of random sum complexes, {\it J. Topol. Anal.},
{\bf 12}(2020) 989–-1002.

\bibitem{Bjorner92}
A. Bj\"{o}rner, The homology and shellability of matroids and geometric lattices. Matroid applications, 226-–283, Encyclopedia Math. Appl., 40, Cambridge Univ. Press, Cambridge, 1992.

\bibitem{Bollobas98}
B. Bollob\'as, {\it Modern Graph Theory}, Graduate Texts in
Mathematics, Springer Verlag, New York, 1998.

\bibitem{DR02}
A. Duval and V. Reiner, Shifted simplicial complexes are Laplacian integral,
{\it Trans. Amer. Math. Soc.}, {\bf 354}(2002) 4313-–4344.

\bibitem{Fulton-Harris91}
W. Fulton and J. Harris, Representation theory. A first course. Graduate Texts in Mathematics, 129.
Readings in Mathematics. Springer-Verlag, New York, 1991.
	
\bibitem{Garland73}
H. Garland, $p$-adic curvature and the cohomology of discrete
subgroups of $p$-adic groups, {\it Annals of Math.}, {\bf 97}(1973) 375--423.	

\bibitem{HRW98}
J. Herzog, V. Reiner and V. Welker, The Koszul property in affine semigroup rings,
{\it Pacific J. Math.}, {\bf 186}(1998) 39-–65.


\bibitem{Landau-Russell04}
Z. Landau and A. Russell, Random Cayley graphs are expanders: a simple proof of the Alon-Roichman theorem,
{\it Electron. J. Combin.}, {\bf 11}(2004), no. 1, Research Paper 62, 6 pp.

\bibitem{LMR10}
N. Linial, R. Meshulam and M. Rosenthal, Sum complexes -
a new family of hypertrees,
{\it Discrete Comput. Geom.}, {\bf 44}(2010) 622--636.

\bibitem{Loh-Schulman04}
P.-S. Loh and L. J. Schulman, Improved expansion of random Cayley graphs, {\it Discrete Math. Theor. Comput. Sci.},
{\bf 6}(2004) 523–-528.

\bibitem{Lubotzky18}
A. Lubotzky, High dimensional expanders. Proceedings of the International Congress of Mathematicians—Rio de Janeiro 2018. Vol. I. Plenary lectures, 705–730, World Sci. Publ., Hackensack, NJ, 2018.

\bibitem{LPS88}
A. Lubotzky, R. Phillips and P. Sarnak, Ramanujan graphs, {\it Combinatorica},
{\bf 8}(1988) 261-–277.

\bibitem{M14}
R. Meshulam, Uncertainty principles and sum complexes, {\it Journal of Algebraic Combinatorics}, {\bf 40}(2014) 887--902.

\bibitem{Tropp12}
J. A. Tropp, User-friendly tail bounds for sums of random matrices, {\it Found. Comput. Math.}, {\bf 12}(2012) 389–-434.

\bibitem{ZZ93}
G.M. Ziegler and R. \v{Z}ivaljevi\'{c}, Homotopy types of subspace arrangements via diagrams
of spaces, {\it Math. Ann.}, {\bf 295}(1993) 527--548.



\end{thebibliography}
\end{document}